\documentclass[12pt,dvips]{amsart}
\usepackage{euler, amsfonts, amssymb, latexsym, epsfig,epic}

\newtheorem{Theorem}{Theorem} 
\newtheorem{Proposition}{Proposition} 
\newtheorem{Lemma}{Lemma}
\newtheorem{Fact}{Fact}

\newtheorem*{Corollary}{Corollary}
 
\newtheorem*{Theorem*}{Theorem}
\theoremstyle{remark}
\newtheorem{Example}{Example}

\newcommand\onto{\mathop{\to}}
\newcommand\into{\operatorname*{\hookrightarrow}}

\newcommand\CP{{\mathbb C \mathbb P}}

\newcommand\complexes{{\mathbb C}}
\newcommand\integers{{\mathbb Z}}

\theoremstyle{plain}

\newcommand\dfn{\bf} % maybe should be \em

\newcommand\GLn{{{GL_n(\complexes)}}}

\begin{document}
\pagestyle{plain}

\title{Descent-cycling in Schubert calculus}
\author{Allen Knutson}
\thanks{This research was partially conducted for 
  the Clay Mathematics Institute.}
\email{allenk@math.berkeley.edu}
\date{\today}

\begin{abstract}
  We prove two lemmata about Schubert calculus on generalized flag
  manifolds $G/B$, and in the case of the ordinary flag manifold
  $GL_n/B$ we interpret them combinatorially in terms of descents, and
  geometrically in terms of missing subspaces.  One of them gives a
  symmetry of Schubert calculus that we christen {\em descent-cycling}.  
  Computer experiment shows that these lemmata suffice to determine
  all of $GL_n$ Schubert calculus through $n=5$, and $99.97\%+$ at
  $n=6$. We use them to give a quick proof of Monk's rule. The lemmata
  also hold in equivariant (``double'') Schubert calculus for
  Kac-Moody groups $G$.
%, as studied in \cite{KK}.
\end{abstract}

\maketitle

\section{Background on Schubert problems}\label{sec:intro}

Fix a pinning for a complex reductive Lie group $G$: a Borel subgroup $B$, 
an opposed Borel subgroup $B_-$, a Cartan subgroup $T = B \cap B_-$,
the Weyl group $W = N(T)/T$, and $R$ the Coxeter generators of $W$. 
There is a famous basis (as a free abelian group) for the cohomology
of $G/B$ given by the Poincar\'e duals of the closures of the $B_-$
orbits on $G/B$; these are the {\dfn Schubert classes} 
$S_w := [\overline{B_- w B/B}], w\in W$, and are indexed by the Weyl group. 

(In this introduction we will only consider ordinary cohomology and the
case of finite-dimensional $G$. However, since the Schubert cycles
$\overline{B_- w B/B}$ are $T$-invariant, they define elements not
only of ordinary but of $T$-equivariant cohomology of $G/B$%, 
%(and we will use the same notation $S_w$ without fear of confusion)
, and our results hold in that case also.
In addition, our main arguments apply to the case of Kac-Moody $G$.
Our references for equivariant cohomology of (possibly infinite-dimensional)
$G/B$ are \cite{G,KK}.)

The degree of the cohomology class $S_w$ is twice $l(w)$, 
the {\em length} of $w$ (as a minimal product of Coxeter generators from $R$).
Define a {\dfn Schubert problem} as a triple $(u,v,w) \in W^3$
such that $l(u)+l(v)+l(w) = \dim_\complexes G/B$.
In this case we can consider the {\dfn symmetric Schubert numbers}
$$ c_{uvw} := \int_{G/B} S_u S_v S_w $$
which count the number of points in the intersection 
of three generic translates of Schubert cycles.
Since this intersection is transverse (by a standard appeal to Kleiman's
transversality theorem), and is of three complex subvarieties, the
points are all counted with sign $+1$ and therefore the number is
nonnegative. It is a famous open problem to compute this number
combinatorially; the analogous problem for $G/P$ where $G=GL_n$ and
$P$ is a maximal parabolic was solved first by the Littlewood-Richardson rule
(or see \cite{KT}).

Recall the Bruhat order on $W$ (due to Chevalley): $v>w$ if 
$v \in \overline{B_- w B/B}$. With this we can state our two lemmata:

\begin{Lemma}\label{lem:dctriv}
  Let $(u,v,w)$ be a Schubert problem, and $s$ a Coxeter generator.
  If $us > u$, $vs > v$, and $ws > w$, then $c_{uvw} = 0$.
\end{Lemma}

\begin{Lemma}\label{lem:dc}
  Let $(u,v,w) \in W^3$ be a triple with 
  $l(u)+l(v)+l(w) = \dim_\complexes G/B -1$, and $s$ a Coxeter generator.
  If $us > u$, $vs > v$, and $ws > w$, then 
  $$ c_{us,v,w} = c_{u,vs,w} = c_{u,v,ws}. $$
\end{Lemma}

In the case $G=GL_n(\complexes), W=S_n,$ and $s$ is the transposition
$i \leftrightarrow i+1$, the statement $us>u$ says that $u(i)<u(i+1)$;
one says that $u$ {\dfn ascends} in the $i$th place. Otherwise if
$u(i)>u(i+1)$ one says that $u$ {\dfn descends} in the $i$th place,
or that it has a descent there. For this reason we christen the symmetry
of lemma \ref{lem:dc}
{\dfn descent-cycling}, and call these three problems {\dfn dc-equivalent}.
Extending this relation by transitivity, we get a very powerful notion
of equivalence for solving Schubert problems; in particular many Schubert
problems are dc-equivalent to ones that fall to lemma \ref{lem:dctriv},
ones which we call {\dfn dc-trivial}.

Define a {\dfn Grassmannian Schubert problem} $(u,v,w)$ to be one
in which $u,v$ each have only one descent, and in the same place,
so named because the relevant integral can be performed on a Grassmannian;
these Schubert problems are well-understood thanks to 
Littlewood-Richardson and other positive rules for their computation.
It is worth pointing out that descent-cycling {\em cannot} be formulated
in the context of Grassmannian problems alone;
\begin{enumerate}
\item descent-cycling a (nontrivial) Grassmannian Schubert problem
  always produces a non-Grassmannian Schubert problem;
\item Grassmannian problems from different Grassmannians (the single
  descent in different places) can be dc-equivalent.
\end{enumerate}

In section \ref{sec:graph} we define a graph whose vertices are Schubert
problems and edges come from descent-cycling; by computer we were able
to determine much about the structure of this graph in small examples.
This we believe to be the main point of interest in the paper -- that two
such simple lemmata suffice to determine so many Schubert numbers.

It is our hope that this symmetry might help guide the search for a 
combinatorial formula for Schubert calculus; a rule generalizing
Littlewood-Richardson (the case that $\pi,\rho$ each have only one
descent, and in the same place) and manifestly invariant under
descent-cycling would have very strong evidence for it.\footnote{%
I circulated a preprint a year ago entitled ``A conjectural rule
for $\GLn$ Schubert calculus'' in terms of {\em puzzles}, generalizing
a Grassmannian theorem from \cite{KT}. Alas, the rule 
conjectured there is not invariant under descent-cycling.}

In section \ref{sec:cohom} we give the nearly-trivial proofs of the
two lemmata, using standard properties of the (equivariant) BGG operators.
We do this in terms of ``Schubert structure constants'' rather than
symmetric Schubert numbers, which seems to be more appropriate for 
equivariant cohomology, and also gives results in the case of $G$
a Kac-Moody group.

In the $\GLn$ case, there is an intuitive geometrical 
interpretation in terms of ``reconstructing forgotten subspaces'';
with this we can also say something about finding the actual flags 
in the intersection in synthetic-geometry terms, which we do 
in section \ref{sec:geometry}.

In section \ref{sec:Monk} we prove Monk's rule via descent-cycling,
to give an example of an interesting Schubert problem that falls 
to these techniques. It would be interesting to see if other known
cases of $c_{uvw} = 0,1$ (such as the Pieri rule \cite{S}) are
consequences of descent-cycling.

We are thankful to Mark Haiman and Peter Magyar for useful comments.

\section{The Schubert problems graph, 
  and its structure for small $\GLn$}\label{sec:graph}

%Before doing so, we wish to comment on the effectiveness of descent-cycling
%for actually doing Schubert calculus. 
Let $\Gamma_n$ be the graph whose vertices are Schubert problems
for $\GLn$, with edges between two Schubert problems that are related by
cycling a single descent. Then the descent-cycling lemma \ref{lem:dc}
says that the symmetric Schubert number is constant on connected 
components of this graph.%
\footnote{Probably
a graph structure is not quite the right one to use for this, since the
natural concept of ``edge'' here connects three, not two, vertices.}
Recall that we define two Schubert problems to be dc-equivalent if they
are in the same connected component, i.e., if one can be transformed
into the other by a sequence of descent-cyclings. Also, we call a Schubert
problem dc-trivial if it falls to lemma \ref{lem:dctriv}, i.e.
for some $(i,i+1)$ it has three ascents.

\begin{Example}
We write a vertical bar to point out the descents, and a horizontal bar
indicating to where we intend to cycle a descent. 
\begin{figure}[htbp]
  \begin{center}
    \epsfig{file=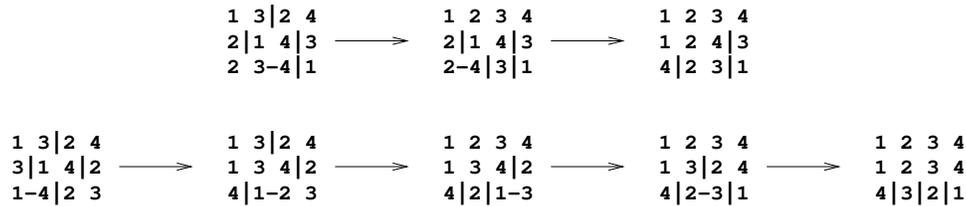,width=5in}
    \caption{Using descent-cycling to simplify a Schubert problem,
      moving the vertical bars (marking descents) onto the horizontal ones.}
    \label{fig:ex1}
  \end{center}
\end{figure}
In the first line of figure \ref{fig:ex1} we descent-cycle our way to
a dc-trivial problem; this shows $c_{1324,2143,2341}=0$. In the second we 
show $(1324,3142,1423)$ dc-equivalent to $(1234,1234,4321)$, so
$c_{1324,3142,1423} = c_{1234,1234,4321}$, which is in turn easily
seen to be $1$. The reader may enjoy studying hands-on the properties of
descent-cycling via the descent-cycling Java applet at

\centerline{
{\tt http://www.math.berkeley.edu/\~{}allenk/java/DCApplet.html}}

\end{Example}

We established the following Facts by brute-force computation.
%the source code for these computations  
%is available at {\tt http://\~{}somewhere}.

\begin{Fact}
There are 35 Schubert problems for $GL_3(\complexes)$, of which 21 are not
dc-trivial.  All 21 live in the same connected component of $\Gamma_3$,
which is pictured in figure \ref{fig:3nontriv}.
\begin{figure}[htbp]
  \begin{center}
    \epsfig{file=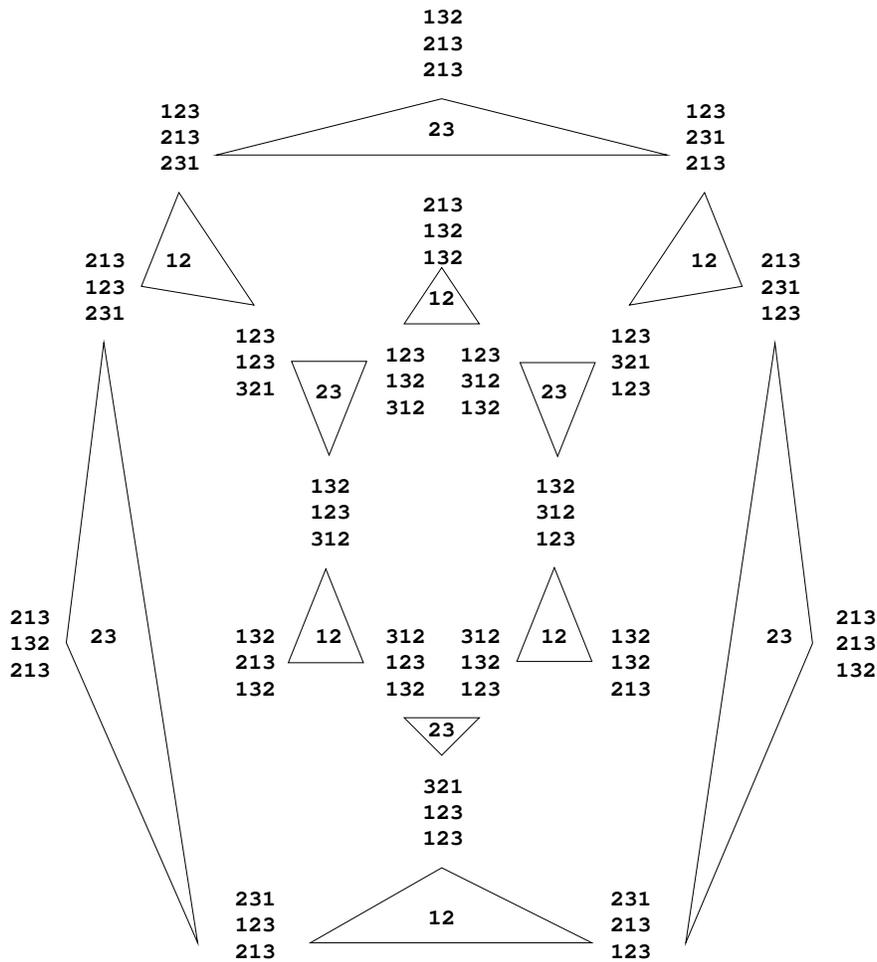,height=5in}
    \caption{The non-dc-trivial component of $\Gamma_3$, drawn to
        make its $S_3$ symmetry manifest. The edges,
        which always come in sets of three, are drawn as triangles
        and labeled with the column where descents are being
        cycled. Note that not all vertices have degree 4; one cannot
        descent-cycle in a column which has {\em two} descents.}
    \label{fig:3nontriv}
  \end{center}
\end{figure}
\end{Fact}

\begin{Fact}
Let $n\leq 5$, and $(\pi,\rho,\sigma)$ be a Schubert problem in dimension $n$.
Then the symmetric Schubert number $c_{\pi\rho\sigma}$ equals zero if 
and only if $(\pi,\rho,\sigma)$ is dc-equivalent to a dc-trivial problem.
Otherwise, $(\pi,\rho,\sigma)$ is dc-equivalent to $(id,id,w_0)$
(where $id$ denotes the identity permutation and $w_0$ the long word),
and therefore the symmetric Schubert number $c_{\pi\rho\sigma}$ equals one.
Put another way, there is exactly one non-dc-trivial component 
of the Schubert problem graph for each $n\leq 5$.
\end{Fact}

In particular, the two lemmata (and the trivial calculation
$c_{id,id,w_0}=1$) suffice to completely determine Schubert calculus 
for $\GLn$ through $n=5$. We know a priori that this connectedness
cannot continue at $n=6$, because
the nonzero symmetric Schubert numbers are sometimes $2$. 
(All symmetric Schubert numbers in this paper were computed with
the Maple package \cite{ACE}.)

\begin{Fact}
The graph $\Gamma_6$ has 8,881,334 vertices, of which 
all but 2,351,475 are dc-trivial. 
Throwing out the components with dc-trivial vertices we are left with
145 components comprising 411,582 vertices.
The lion's share of those vertices, 409,023, are dc-equivalent to the
easy case $(id,id,w_0)$, leaving 2559 cases (less than $0.03\%$) not
succumbing to dc-equivalence/dc-triviality arguments.

Of the remaining 144 components, only one contains a Grassmannian 
Schubert problem, so in some sense the Littlewood-Richardson rule 
doesn't help much. (In fact this is the only one with any
Grassmannian permutations, so a Schubert-times-Schur rule wouldn't
help much either.)

Exactly one of these components has intersection number zero
(despite containing no dc-trivial Schubert problems); 
one element of it is $(231645,231645,326154)$.
There are also 48 components of size one, i.e. Schubert problems that
admit no descent-cycling whatever; one example is $(214365,154326,321654)$.
\end{Fact}

These computations were done in C and took 2.5 minutes on a Pentium 300.
The limiting factor was that they just barely fit in 64 megabytes of RAM,
putting the $n=7$ case (which is roughly $7^3=343$ times bigger) out of 
reach without new ideas. 
%(In any case all our hopeful conjectures,
%like ``$c=0$ implies dc-equivalence to a dc-trivial problem,''
%were already struck down at $n=6$.)

It seems likely that as $n$ increases, the fraction of $\GLn$ Schubert
problems having no place with three ascents (and so falling to lemma
\ref{lem:dctriv} alone) goes to 0. We did not pursue this.

It was very tempting to believe that a vanishing Schubert number 
could always be ``blamed'' on dc-equivalence to a dc-trivial problem,
and it was very sad to find the lone component in $\Gamma_6$ that belies this.
Mark Haiman pointed out a ``stabilization'' map
$\Gamma_{n-1} \into \Gamma_n$ taking $(u,v,w) \mapsto (u n, v n, 1 (w+1))$ 
(where $w+1$ means to add 1 to all elements of $w$). (For example,
$(2143,1243,3214) \mapsto (21435,12435,14325)$.)
Question: can distinct components become connected under this map
(possibly connecting the rogue component in 
$\Gamma_6$ to a dc-trivial problem)?

\section{Proofs of the lemmata}\label{sec:cohom}

The statements in this section are slightly different from those in
the introduction, in that they are phrased in terms of 
{\em structure constants} $c_{uv}^w$, rather than symmetric Schubert numbers
$c_{uvw}$. We first remind the reader of the partial relation between these
and explain why we switch to the less-symmetric formulation.

In ordinary cohomology of $G/B$ ($G$ finite dimensional), we have the
Poincar\'e-pairing duality relation
$$ \int_{G/B} S_u S_v = \delta_{u,w_0 v} $$
where $w_0$ is the long element of the Weyl group, and $\delta$ is
the Kronecker delta. 
(Note that the duality discussed here is in the sense of dual bases,
and not Poincar\'e duality!)

One way to see this is to realize the class $S_v$
not by the Schubert cycle $\overline{B_- v B}$ but the opposite Schubert
cycle $\overline{B w_0 v B} = w_0 \overline{B_- v B}$.
Since $w_0$ is connected to the identity in $G$, these two cycles define
the same element of cohomology.

From this we derive that if
$$ S_u S_v = \sum_w c_{uv}^w\, S_w, \qquad \hbox{ then } \quad
        c_{uv}^w = c_{u,v,w_0 w}. $$
So in ordinary cohomology, one can work instead with these
{\dfn Schubert structure constants}, although our
results from the first section are prettier to state symmetrically.
However, in $T$-equivariant cohomology the dual basis to the Schubert basis
is {\em not} once again the Schubert basis (essentially because $w_0$
is not connected to the identity {\em through $T$-invariant maps} of $G/B$)
and these two concepts part ways.

In \cite{G} a certain positivity result was proven for the equivariant 
$c_{uv}^w \in H^*_T$ (which must be carefully stated, insofar as these
are polynomials not numbers). This implies a much weaker positivity result for
the $c_{uvw}$, and so it seems more interesting to prove results about 
the structure constants.

Also, in the case of $G$ an infinite-dimensional Kac-Moody group, one
cannot so blithely do integrals on $G/B$, and the $c_{uv}^w$ are the
only concept that makes sense.  This concludes the advertisement for
Schubert structure constants over symmetric Schubert numbers.\footnote{%
In fact almost all work on Littlewood-Richardson is in terms of the
structure constants; see \cite{KT} for a discussion of this.}

Note that the condition we gave in section \ref{sec:intro} for a ``Schubert
problem'' corresponds to $l(w)=l(u)+l(v)$, which seems a reasonable
thing to ask since cohomology is a graded ring. But we will not insist on
this in what follows because, in equivariant cohomology, the structure 
constants can be nonzero even if one only has $l(w) \leq l(u)+l(v)$.
(We only imposed this condition before to keep the graph of Schubert problems
a reasonable size.)

Let $s\in R$ be a simple reflection, and $P_s = \overline{BsB} \leq G$ be 
the corresponding minimal parabolic. Let $p_s : G/B \onto G/P_s$ be the
corresponding projection (which is $G$-equivariant and therefore 
$T$-equivariant); composing pushforward with pullback 
defines a degree $-2$ endomorphism $\partial_s$ on $H^*_T(G/B)$ 
(first introduced in \cite{BGG} and \cite{D}, stated nonequivariantly,
although equivariant K-theory is implicit in \cite{D}). 

We refer to \cite{KK} for the four properties we need of
these BGG operators $\partial_s$:
\begin{enumerate}
\item If $ws > w$, then $\partial_s S_w = 0$.
\item If $ws < w$, then $\partial_s S_w = S_{ws}$.
\item With respect to a certain natural action of $W$ on $H^*_T(G/B)$, 
$\partial_s$ is a twisted derivation:
$$ \partial_s (\alpha \beta) 
= \alpha \, \partial_s(\beta) + \partial_s(\alpha) \, (s\cdot \beta). $$
(We won't need to understand this action of $W$.)
\item If $s_1 s_2 \ldots s_l = r_1 r_2 \ldots r_l$ are two 
reduced expressions for a Weyl group element $w$, then
$\partial_{s_1} \ldots \partial_{s_l} = \partial_{r_1} \ldots \partial_{r_l}$,
and so we have a well-defined operator $\partial_w$.
\end{enumerate}

Since the proofs of both lemmata have much in common, we gather them into a
single proposition.

\begin{Proposition}\label{prop:dcboth}
  Let $(u,v,w)\in W^3$, and $s$ a simple reflection, 
  such that $us>u$ but $ws<w$.
  \begin{itemize}
  \item If $vs>v$, then $c_{uv}^w = 0$.\label{i:dctriv}
  \item If $vs<v$, then $c_{uv}^w = c_{u,vs}^{ws}$.\label{i:dc}
  \end{itemize}
\end{Proposition}

\begin{proof}
  The main formula we need is
  $$ c_{xy}^z = \hbox{coefficient of $S_1$ in } \partial_{z} (S_x S_y) $$
  which follows easily from the properties stated of the BGG operators.

  Since $ws<w$, we have
  $$ \partial_w (S_u S_v) 
  = \partial_{ws} \partial_s (S_u S_v) 
  = \partial_{ws} ( S_u \partial_s (S_v) + \partial(S_u) (s\cdot S_v) )
  = \partial_{ws} ( S_u \partial_s (S_v) ) $$
  this last because by the $us>u$ assumption, $\partial_s$ annihilates $S_u$.
  
  If $vs>v$, then $\partial_s$ annihilates $S_v$ too, and the RHS is zero.
  Combining that with the formula for $c_{xy}^z$ gives the first result.

  If $vs<v$, then the RHS is $\partial_{ws} (S_u S_{vs})$, and two 
  applications of the formula for $c_{xy}^z$ give the second result.
\end{proof}

The conditions on $ws$ versus $w$ in this proposition are backwards from 
how they were in lemmata \ref{lem:dctriv} and \ref{lem:dc}; that's
because of the multiplication by $w_0$ in comparing Schubert structure
constants with symmetric Schubert numbers. With this in mind the
two lemmata follow.

\section{A geometrical interpretation}\label{sec:geometry}

For $w\in W$, let $D_w := \overline{G \cdot (w B,B)} \subseteq (G/B)^2$.
Given a simple reflection $s\in R$, let $P_s$ again be the corresponding
minimal parabolic $\overline{BsB}$, and consider the composite map
$D_w \into (G/B)^2 \onto G/P_s \times G/B$.

\begin{Lemma}
  Let $w\in W, s\in R, P = \overline{BsB}$. The fibers of %the composite map
  $D_w \into (G/B)^2 \onto G/P_s \times G/B$ are
\begin{itemize}
\item $\CP^1$'s if $ws < w$
\item single points (generically), if $ws > w$.
\end{itemize}
\end{Lemma}

\begin{proof}
We reduce to the well-studied case (see \cite{D}) of a single flag manifold.
Let $X := G/B \times \{B/B\}$; since $G \cdot X = G/B \times G/B$ it 
suffices to consider the map $D_w \cap X \to G/P_s \times \{B/B\}$.
And $D_w \cap X = \overline{BwB/B} \times \{B/B\}$, 
so (omitting the $\{B/B\}$) we're studying the fibers of the composite
$\overline{BwB/B} \to G/B \onto G/P_s$, as already done in \cite{D}.
\end{proof}

In the case of $G = \GLn$, % and denoting $G/B$ by $Fl_n$ for flag manifold, 
$D_w$ is the variety of pairs of flags $(F,G)$ in $\complexes^n$ such that
``$F$ is $w$-close or closer to $G$''. In this case, the generators
$R$ correspond 1:1 to the subspaces in a flag (other than the zero subspace
and the whole space), and the map $G/B \to G/P_s$ corresponds to 
``forgetting'' the subspace. Then we can interpret the lemma 
in very familiar terms:

\begin{Corollary}
  Let $w\in S_n$, $i \in 2\ldots n-1$. Let $F,G$ be two flags in
  $\complexes^n$ such that $F$ is $w$-close or closer to $G$.  Let $F'$ be 
  the partial flag obtained by forgetting $F$'s $i$-dimensional
  subspace. Can we reconstruct $F$ knowing only $F'$, $G$, and $w$?
\begin{itemize}
\item If $w$ ascends at $(i,i+1)$, there is no hope -- 
  any $i$-dimensional space in between $F_{i-1}$ and $F_{i+1}$ will do.
\item If $w$ descends at $(i,i+1)$, then (for a Zariski-open set of
  such $G$) the subspace $F_i$ is uniquely determined.
\end{itemize}
\end{Corollary}

Another way to interpret this is that if $w$ does not descend at $(i,i+1)$,
then $G$ ``does not care'' what $F_i$ is used (to get $F$ $w$-close 
to $G$). Conversely, if $w$ does descend there, 
then $G$ ``usually insists'' on a particular $F_i$, when presented with
the rest of $F$.

\begin{proof}[Proof of lemma \ref{lem:dctriv}]
Let $\{F\}$ be the set of flags in relative position $u$ to $A$,
$v$ to $B$, and $w$ to $C$ where $A,B,C$ are three flags
in generic relative position. Then by codimension count (and the
usual appeal to Kleiman's transversality theorem) the set $\{F\}$ is finite.
However, since none of $A,B,C$ care what $F_i$ is (since by assumption
none of them have a descent at $(i,i+1)$), the set $\{F\}$ is a union
of $\CP^1$'s. These two facts are only compatible if $\{F\}$ is empty.
\end{proof}

\begin{proof}[Proof of lemma \ref{lem:dc}]
Let $\{F\}$ be the set of flags in relative position $u$ to $A$,
$v$ to $B$, and $w$ to $C$ where $A,B,C$ are three flags
in generic relative position. Then as in the previous proof, the set
$\{F\}$ is a union of $\CP^1$'s, reflecting the ambiguity in $F_i$. 
If we change one of $u,v,w$ to have a descent at $(i,i+1)$, each
of these $\CP^1$'s is cut down to a single point. But it doesn't matter
which of $u,v,w$ gets this new descent.
\end{proof}

This geometric description of descent-cycling suggests that additional
symmetries may come from forgetting multiple subspaces at a time.
It appears, though, that all of these are implied by the single-subspace case.

One application of this geometric description is to actually locate
the flag satisfying the desired intersection conditions, in the case
that $(\pi,\rho,\sigma)$ is dc-equivalent to the easy case $(id,id,w_0)$. 
We illustrate this in the case of the Schubert problem $(132,213,213)$,
which we can descent-cycle to $(123,213,231)$, and from 
there to $(123,123,321)$. 
Working from the end, the unique $F$ satisfying $(123,123,321)$ 
is $F_1=C_1, F_2=C_2$.
When we cycle the descent in the $(1,2)$ column, we have to replace
$F_1=C_1$ by $F_1=B_2\cap C_2$. Then when we cycle the descent in the $(2,3)$
column we have to replace $F_2=C_2$ by $F_2 = A_1\oplus (B_2\cap C_2)$.

There is an alternate way to prove the vanishing condition in
proposition \ref{prop:dcboth} cohomologically, involving the projection
$G/B \to G/P_s$. A Schubert class $S_u$ is in the image of the
pullback of $H^*_T(G/P_s)$ if and only if $us > u$. Since this pullback
is a ring homomorphism, the product of two pulled-back classes is also in
this image, and cannot involve any $S_w$ with $ws < w$.

\section{Monk's rule}\label{sec:Monk}

Monk's rule \cite{M} is concerned with the case of $\GLn$ Schubert problems
in which $\rho$ is a simple reflection $s_i = (i\leftrightarrow i+1)$. 
%(It is the $GL_n$ case of a rule discovered at the same time by
%Chevalley \cite{C} for all $G$, about which we will have nothing to say.)

\begin{Theorem*}[Monk's rule]
Let $\sigma w_0$ cover $\pi$ in the Bruhat order; i.e. $\sigma$ is $\pi$
with each number $j$ replaced by $(n+1)-j$, and two numbers inverted
in $\pi w_0$ have been put back in correct order, 
decreasing the number of inversions by exactly one. Then 
$c_{\pi,s_i,\sigma}=1$ if the numbers switched straddled the position
between $i$ and $i+1$, whereas
$c_{\pi,s_i,\sigma}=0$ if the numbers switched were both physically 
on one side of $(i,i+1)$. 
\end{Theorem*}

(Some may object that Monk's rule says more -- that $c_{\pi,s_i,\sigma}=0$
unless $\sigma w_0$ covers $\pi$ -- but we prefer to see this as a more
general property %(see \cite{}) 
of symmetric Schubert numbers, that
if $\pi,\rho$ are not less than $\sigma w_0$ in the Bruhat order,
then $c_{\pi\rho\sigma}=0$.) 

For example, let $\pi=34152, i=2$. Then $\pi w_0 = 32514$, which
covers $23514$, $31524$, $32154$, $32415$. But only $31524$ involves
switching a number in the first $2$ places with a number in the last
$5-2$ places.  So $c_{34152,s_i,31524}=1$, but $c_{34152,s_i,23514} =
c_{34152,s_i,32154} = c_{34152,s_i,32415} = 0$.

\begin{Theorem}\label{thm:Monk}
Let $f : V(\Gamma_n) \to \integers$ be a functional on the set of 
Schubert problems. If $f$ satisfies the properties
\begin{enumerate}
%\item $f(\pi,\rho,\sigma)=0$ unless $\pi \geq \sigma w_0$ in the Bruhat order
\item $f$ is invariant under descent-cycling (i.e. is constant on components)
\item $f=0$ on dc-trivial Schubert problems
\item $f(id,id,w_0)=1$
\end{enumerate}
then $f$ obeys Monk's rule, i.e. $f(\pi,s_i,\sigma)=0$ or $1$ according to
the criterion of Monk's rule.
\end{Theorem}

We first prove a lemma:

\begin{Lemma}\label{lem:Monk}
Let $f$ satisfy the conditions of theorem \ref{thm:Monk}, 
and $\pi,\sigma \in S_n$ such that $l(\pi)+l(\sigma)={n\choose 2}$.
Then if $\pi=\sigma w_0$, we have $f(\pi,id,\sigma)=1$, and otherwise
$f(\pi,id,\sigma)=0$.
\end{Lemma}

\begin{proof}
Since the second argument has no descents, any place $(i,i+1)$
that $\pi$ has a descent and $\sigma$ does not 
gives us an opportunity to cycle a descent from
the first argument to the third, replacing $\pi\mapsto s_i \pi,
\sigma \mapsto s_i \sigma$. This modification keeps the sum of the
lengths $= {n\choose 2}$ and neither causes nor breaks the condition
$\pi=\sigma w_0$. So we can reduce to the case that any descent in
$\pi$ occurs at a descent of $\sigma$.

If $\pi=\sigma w_0$: then each ascent in $\sigma$ occurs at a descent
of $\pi$. By our reduction above, this means that $\sigma$ has no ascents.
So we're looking at $f(id,id,w_0)$ which by assumption is $1$.

Conversely if $f(\pi,id,\sigma)\neq 0$: then no column $(i,i+1)$ has
three ascents (the Schubert problem $(\pi,id,\sigma)$ is not
dc-trivial).  By our reduction, this means that $\sigma$ has no
ascents. So $\sigma=w_0$.  By the assumption on the total length,
$\pi=id$, so $\pi=\sigma w_0$ as desired.
\end{proof}

\begin{proof}[Proof of theorem \ref{thm:Monk}]
Let $\sigma$ be $\pi w_0$ with the numbers in the $j$th and $k$th positions
switched, decreasing the number of inversions by exactly one (and so that
$j<k, \sigma(j)<\sigma(k)$). In particular every number in $\sigma$
physically between the $j$th and $k$th positions is {\em not}
numerically between $\sigma(j)$ and $\sigma(k)$.
We want to show that $f(\pi,s_i,\sigma)=0$
unless $j\leq i < k$, in which case $f(\pi,s_i,\sigma)=1$.

First, we treat the case $k=j+1$.
If $j=i, k=i+1$, then neither $\pi$ nor $\sigma$ has a descent at $(i,i+1)$.
So we can cycle the descent from the second argument of $f$ into
the third, making them $(\pi,id,\pi w_0)$. Now lemma \ref{lem:Monk} 
tells us that this $1$.

Whereas if $j,k\leq i$ or $j,k\geq i+1$ (still in the case $k=j+1$), 
then none of $\pi$, $id$, or $\sigma$ have a descent at $(j,k)$,
and therefore $f$ vanishes as it's supposed to.

Now take the case $k>j+1$. Then since $\sigma$ has only one fewer inversion
than $\pi w_0$, $\sigma$ must have the same descent-pattern as $\pi w_0$.
Now we reduce (much as in lemma \ref{lem:Monk}) by cycling descents
between the first and third arguments, in order to move the positions
$j$ and $k$ closer together. 

We can do this descent-cycling in the $(j,j+1)$ column as long as
$j\neq i$, and the $(k-1,k)$ column as long as $k\neq i+1$. If $j,k$ are
both on the same side of the $(i,i+1)$ divide, they can be brought next
to each other (by e.g. just moving one of them). If $j,k$ are on 
opposite sides of the divide, we can at least get $j$ up to $i$, and
$k$ down to $i+1$. Either way we reduce to the $k=j+1$ case and therefore
get the same answer as Monk's rule.
\end{proof}

In particular, this gives an explicit sequence of descent-cyclings to
turn a Monk's rule problem into $(id,id,w_0)$. So in principle one can
reverse the steps and construct the flag in the intersection of these
three Schubert varieties, as an expression in the lattice of subspaces.

There are other special cases known for symmetric Schubert numbers
where the answer is $0$ or $1$, mostly notably the Pieri rule 
(see \cite{R,S}); it would be interesting to see if they too are 
consequences of descent-cycling. Probably the best version of this
would be a ``descent-cycling normal form'' for Schubert problems,
and an effective way to test whether a Schubert problem is dc-equivalent
to $(id,id,w_0)$. 

\section{Syntheticity vs. $c=1$ questions}

Recall that given a Schubert problem $P = (\pi,\rho,\sigma)$,
and three generically situated flags $A,B,C$, 
one can think of the symmetric Schubert number $c_{\pi\rho\sigma}$
as the number of flags $F$ such that $F$ is $\pi$-close to $A$,
$\rho$-close to $B$, and $\sigma$-close to $C$.

The Schubert problem $(id,id,w_0)$ is then easily seen to have
symmetric Schubert number one; to be $id$-close to $A$ or $B$ is no
condition at all, and to be $w_0$-close to $C$ requires $F_i = C_i$
for all $i=1,\ldots,n$.

Consider the following four statements one might make about a
Schubert problem $P = (\pi,\rho,\sigma)$:

\begin{itemize}
\item {\dfn dc-easiness:} $P$ is dc-equivalent to the Schubert problem
  $(id,id,w_0)$
\item {\dfn partial syntheticity:} there is a flag in the free modular
  lattice on three flags $A,B,C$ satisfying $P$
\item {\dfn full syntheticity:} every flag satisfying $P$ is in the free
  modular lattice on three flags $A,B,C$
\item {\dfn $c=1$:} the symmetric Schubert number $c_{\pi\rho\sigma} = 1$.
\end{itemize}

So $P$ dc-easy implies $P$ fully synthetic and $c_P = 1$. The other
possible implications seem to be unknown.

{\em Question.} Does $c=1$ imply partial (and thus full) syntheticity?
This would seem to be a Galois theory argument, with ``synthetic'' the
analogue of ``rational.''

{\em Question.} Does $P$ partially synthetic imply $P$ stably dc-easy?
If the flag $F$ is a synthetic solution to the Schubert problem $P$, 
i.e. $F_i$ is a lattice word in $A,B,C$ for each $i$, 
perhaps there is an algorithm to ``simplify'' the ``worst'' subspace in 
$P$ using descent-cycling, with the only unsimplifiable subspaces being
those in $A,B,C$. In particular, this would say that $P$ partially
synthetic implies $P$ fully synthetic.

\bibliographystyle{alpha}    % it seems this does nothing.

\end{document}